\documentclass[a4paper,10pt,reqno]{amsart}
\usepackage{latexsym}
\usepackage{ae}

\def\disp{\displaystyle}
\newcommand{\pref}[1]{{(\protect\ref{#1})}}

\def\newop#1{\expandafter\def\csname #1\endcsname{\mathop{\rm #1}\nolimits}}
\newop{nlabel}
\newop{comp}
\newop{leaves}
\newop{lroot}
\newop{runs}
\newop{lpath}
\newop{rpath}
\newop{lrmax}
\newop{lrmin}
\newop{rlmax}
\newop{rlmin}
\newop{nonleaves}
\newop{zstat}
\newop{stem}
\newop{lir}
\newop{ldr}

\newcommand{\btree}{\mbox{\ensuremath{\beta(1,0)}-tree}}
\newcommand{\btrees}{\btree s}

\catcode`\@=11
\def\section{\@startsection{section}{1}%
 \z@{.7\linespacing\@plus\linespacing}{.5\linespacing}%
 {\normalfont\bfseries\scshape\centering}}
\def\subsection{\@startsection{subsection}{2}%
  \z@{.5\linespacing\@plus\linespacing}{.5\linespacing}%
  {\normalfont\bfseries\scshape}}

\def\subsubsection{\@startsection{subsubsection}{3}%
  \z@{.5\linespacing\@plus.7\linespacing}{-.5em}%
  {\normalfont\itshape}}
\catcode`\@=12
\theoremstyle{plain}

\newtheorem*{theorem*}{Theorem}

\newtheorem*{corollary*}{Corollary}

\newtheorem*{lemma*}{Lemma}

\newtheorem*{proposition*}{Proposition}

\newtheorem*{conjecture*}{Conjecture}

\theoremstyle{definition}

\newtheorem*{definition*}{Definition}

\newtheorem*{example*}{Example}

\newtheorem*{problem*}{Problem}

\theoremstyle{remark}

\newtheorem*{remark*}{Remark}

\def\dd{\makebox[1.1ex]{\rule[.58ex]{.71ex}{.15ex}}}

\def\yab{\ensuremath{[12)}}

\def\ba{\ensuremath{21}}
\def\bxa{\ensuremath{2{\dd}1}}
\def\abc{\ensuremath{123}}
\def\axbc{\ensuremath{1{\dd}23}}

\def\axbxc{\ensuremath{1{\dd}2{\dd}3}}
\def\acb{\ensuremath{132}}
\def\axcb{\ensuremath{1{\dd}32}}
\def\acxb{\ensuremath{13{\dd}2}}
\def\axcxb{\ensuremath{1{\dd}3{\dd}2}}
\def\bac{\ensuremath{213}}
\def\bxac{\ensuremath{2{\dd}13}}

\def\bca{\ensuremath{231}}
\def\bxca{\ensuremath{2{\dd}31}}

\def\bxcxa{\ensuremath{2{\dd}3{\dd}1}}
\def\cab{\ensuremath{312}}

\def\caxb{\ensuremath{31{\dd}2}}

\def\cba{\ensuremath{321}}

\def\cbxa{\ensuremath{32{\dd}1}}

\def\axcxbxd{\ensuremath{1{\dd}3{\dd}2{\dd}4}}

\def\axcdb{\ensuremath{1{\dd}342}}

\def\bxadxc{\ensuremath{2{\dd}14{\dd}3}}

\def\bxaxdxc{\ensuremath{2{\dd}1{\dd}4{\dd}3}}

\def\bxcxda{\ensuremath{2{\dd}3{\dd}41}}

\def\bxcxdxa{\ensuremath{2{\dd}3{\dd}4{\dd}1}}

\def\bxdaxc{\ensuremath{2{\dd}41{\dd}3}}

\def\bxdxaxc{\ensuremath{2{\dd}4{\dd}1{\dd}3}}

\def\cxabxd{\ensuremath{3{\dd}12{\dd}4}}

\def\cxadxb{\ensuremath{3{\dd}14{\dd}2}}

\def\cxaxdxb{\ensuremath{3{\dd}1{\dd}4{\dd}2}}

\def\cxbxda{\ensuremath{3{\dd}2{\dd}41}}

\def\cxdaxb{\ensuremath{3{\dd}41{\dd}2}}

\def\dxabxc{\ensuremath{4{\dd}12{\dd}3}}

\def\dxbac{\ensuremath{4{\dd}213}}

\def\emm#1,{{\em #1}}
\def\ch#1,#2,{#1\choose#2}

\def\steingr{Steingr\'imsson}

\thanks{The work presented here was supported by grant
  no.\ 060005012/3 from the Icelandic Research Fund}

\title{Generalized permutation patterns --- a short survey}

\author[E.  Steingr\'{\i}msson]{Einar Steingr\'{\i}msson}

\begin{document}
\begin{abstract}
An occurrence of a classical pattern $p$ in a permutation $\pi$ is a
subsequence of $\pi$ whose letters are in the same relative order (of
size) as those in~$p$.  In an occurrence of a \emm generalized
pattern,, some letters of that subsequence may be required to be
adjacent in the permutation.  Subsets of permutations characterized by
the avoidance---or the prescribed number of occurrences---of
generalized patterns exhibit connections to an enormous variety of
other combinatorial structures, some of them apparently deep.  We give
a short overview of the state of the art for generalized patterns.

\end{abstract}

\maketitle
\thispagestyle{empty}

\section{Introduction}\label{section-intro}

Patterns in permutations have been studied sporadically, often
implicitly, for over a century, but in the last two decades this area
has grown explosively, with several hundred published papers.  As
seems to be the case with most things in enumerative combinatorics,
some instances of permutation patterns can be found already in
MacMahon's classical book from 1915, {\em Combinatory Analysis}
\cite{macmahon}.  In the seminal paper \emm Restricted permutations,
of Simion and Schmidt~\cite{simion-schmidt} from 1985 the systematic
study of permutation patterns was launched, and it now seems clear
that this field will continue growing for a long time to come, due to
its plethora of problems that range from the easy to the seemingly
impossible, with a rich middle ground of challenging but solvable
problems.  Most important, perhaps, for the future growth of the
subject, is the wealth of connections to other branches of
combinatorics, other fields of mathematics, and to other disciplines
such as computer science and physics.

Whereas an occurrence of a classical pattern $p$ in a permutation
$\pi$ is simply a subsequence of $\pi$ whose letters are in the same
relative order (of size) as those in~$p$, in an occurrence of a \emm
generalized pattern,, some letters of that subsequence may be required
to be adjacent in the permutation.  For example, the classical pattern
$1\dd2\dd3\dd4$ simply corresponds to an increasing subsequence of
length four, whereas an occurrence of the generalized pattern
$1\dd23\dd4$ would require the middle two letters of that sequence to
be adjacent in $\pi$, due to the absence of a dash between 2 and 3.
Thus, the permutation $23145$ contains $1\dd2\dd3\dd4$ but not
$1\dd23\dd4$.  Note that for the classical patterns, our notation
differs from the usual one, since the dashes we have between every
pair of adjacent letters in a classical pattern are usually omitted
when only classical patterns are being considered.


It is well known that the number of permutations of length $n$
avoiding any one classical pattern of length 3 is the $n$-th Catalan
number, which counts a myriad different combinatorial objects.  There
are many other results in this direction, relating pattern avoiding
permutations to various other combinatorial structures, either via
bijections, or by such classes of permutations being equinumerous to
the structures in question without there being a known bijection.
Counting permutations according to the number of occurrences of
generalized patterns one comes up with a vast array of known
sequences, such as the Euler numbers, Stirling numbers of both kinds,
Motzkin numbers, Entringer numbers, Schr\"oder numbers, Fibonacci
numbers, Pell numbers and many more.  Also, one often finds lesser
known sequences that are nevertheless related to known structures,
such as directed animals, planar maps, permutation tableaux, various
kinds of trees and involutions in $S_n$, to name a few.  Thus,
generalized patterns provide a significant addition to the already
sizable flora of classical patterns and their connections to other
combinatorial structures.

In fact, due to their great diversity, the non-classical generalized
patterns are likely to provide richer connections to other
combinatorial structures than the classical ones do.  Supporting this
is the fact that the recently proved Stanley-Wilf conjecture---which
gives a strong bound for the growth rate of the number of permutations
of length $n$ avoiding a classical pattern---does not hold for some
generalized patterns.  

This paper is organized as follows: In Section \ref{sec-def} we
introduce definitions and in Section \ref{sec-lit} we mention implicit
appearances of generalized patterns in the literature. In Sections
\ref{sec-gp3} and \ref{sec-gp4} we survey what is known about the
avoidance of generalized patterns of length three and four,
respectively.  In Section \ref{section-context} we give some examples
where generalized patterns have shown up in very different contexts,
establishing connections to various other combinatorial structures,
some of which seem quite deep.  Section \ref{section-disguise} lists
several instances of so-called barred patterns that turn out to be
equivalent to generalized patterns and Section \ref{sec-asymp} deals
with asymptotics for avoidance of generalized patterns.  Finally, in
Section \ref{sec-further}, we mention some further generalizations of
the generalized patterns.

\section{Some definitions}\label{sec-def}

If a permutation $\pi=a_1a_2\dots a_n$ contains the pattern $\axbc$
then clearly the \emm reverse, of $\pi$, that is $a_na_{n-1}\dots
a_1$, contains the reverse of $\axbc$, which is the pattern $\cbxa$.
Since taking the reverse of a permutation is a bijection on the set of
permutations of length $n$, the number of permutations avoiding a
pattern $p$ equals the number of permutations avoiding the reverse of
$p$.  More generally, the distribution---on the set of permutations of
length $n$---of the number of occurrences of a pattern $p$ equals the
distribution for the reverse of $p$.  The same is true of the
bijection sending a permutation $\pi=a_1a_2\dots a_n$ to its \emm
complement, $\pi^c=b_1b_2\dots b_n$, where $b_i=n+1-a_i$. (When we
take the complement of a pattern we leave the dashes in place, so the
complement of $\axcdb$ is $\dxbac$.)  These two transformations,
together with their compositions, generate a group of order 4 on the
set of patterns, and we say that two patterns belong to the same \emm
symmetry class, if one is transformed into the other by an element of
this group.  As an example, the patterns $\bxca,\bxac,\acxb$ and
$\caxb$ form an entire symmetry class.\footnote{In the case of
classical patterns, taking the inverse of a pattern is well defined
and preserves avoidance, so patterns that are each other's inverses
belong to the same symmetry class.  This is not the case for
generalized patterns.}

Clearly, two patterns in the same symmetry class have the same
properties with respect to the number of permutations avoiding them,
and more generally when it comes to the number of permutations with
$k$ occurrences of a pattern, for any $k$.  However, it does happen
that two patterns not belonging to the same symmetry class have the
same avoidance.  Thus, two patterns are said to belong to the same
\emm Wilf class, if they have the same avoidance, that is, if the
number of permutations of length $n$ that avoid one is the same as
that number for the other.  Clearly, a Wilf class is a union of
symmetry classes (if we define both as equivalence classes).

For example, although the classical patterns of length 3 belong to two
symmetry classes (represented by $\axbxc$ and $\axcxb$, respectively),
they all belong to the same Wilf class, since their avoidance is given
by the Catalan numbers.  For classical patterns, much is known about
Wilf classes for patterns of length up to 7, but a general solution
seems distant.  For the generalized patterns, much less is known.  The
best reference to date is probably \cite{cla-man-enum-bast}.  As an
example, the patterns $\axbc$ and $\axcb$ belong to the same Wilf
class (the avoidance counted by the Bell numbers in both cases) but
not the same symmetry class.

\section{Generalized patterns  in the literature}\label{sec-lit}

Generalized patterns have shown up implicitly in the literature in
various places, and subsets of these have been studied in some
generality.  Namely, Simion and Stanton \cite{simstant} essentially
studied the patterns \bxca, \bxac, \caxb, and \acxb\ and their
relation to a set of orthogonal polynomials generalizing the Laguerre
polynomials, and one of these patterns also played a crucial role in
the proof by Foata and Zeilberger \cite{FZ} that Denert's statistic is
Mahonian.  A permutation statistic is Mahonian if it has the same
distribution---on the set of permutations of length $n$ for each
$n$---as the number of inversions.

Goulden and Jackson give an exponential generating function (EGF) for
the number of permutations avoiding the pattern 123 (no dashes), in
the book {\em Combinatorial Enumeration} \cite[Exercise 5.2.17a,
p. 310]{goulden-jackson}.  The formula is
\begin{equation}\label{eqgj}
\left(\disp\sum_{n\ge0}\frac{x^{3n}}{(3n)!}-
\disp\sum_{n\ge0}\frac{x^{3n+1}}{(3n+1)!}\right)^{-1}
\end{equation}
Although this does not seem to be mentioned in \cite{goulden-jackson},
the obvious generalization holds.  Namely, the EGF for the number of
permutations avoiding the dashless pattern $12\cdots k$ is obtained by
replacing $3$ by $k$ in \pref{eqgj}.  It is pointed out in
\cite[Section 3]{kit-pogp} that this general result can be obtained
through an inclusion/exclusion argument similar to one given in
\cite{kit-pogp}.

The dashless patterns of length 3 also appeared earlier, implicitly,
as the \emm valleys, ($\bac$ and $\cab$), the \emm peaks, ($\acb$ and
$\bca$), the \emm double ascents, $\abc$ and the \emm double descents,
$\cba$ in a permutation, the study of which was pioneered by
Fran\c{c}on and Viennot \cite{fravie}, and which is intimately related
to Flajolet's \cite{flajolet} generation of Motzkin paths by means of
certain continued fractions.  This will be mentioned later, in Section
\ref{section-context}, in connection with related, recent developments.

Also, the alternating permutations, which have been studied for a long
time \cite{andre, andre-memoir}, are permutations that avoid the
patterns $123$ and $321$, (with the additional restriction that the
first two letters of the permutation be in decreasing order; otherwise
they are \emm reverse alternating,.  In fact, this extra restriction
is equivalent to the avoidance of the pattern $\yab$ as defined in
\cite{babstein}).

Generalized patterns were first defined explicitly, in full
generality, in the paper {\em Generalized permutation patterns and a
  classification of the Mahonian statistics} \cite{babstein}, where it
was shown that almost all Mahonian permutation statistics in the
literature at that time (up to a certain bijective correspondence
translating ``excedance based'' statistics to ``descent based''
statistics; see \cite{babstein, mad}) could be written as linear
combinations of generalized patterns.  All but one of these statistics
could be expressed as combinations of patterns of length at most 3.
The odd one out was a statistic defined by Haglund \cite{haglund},
which, after translation by the bijection mentioned above, could be
expressed as a combination of patterns of length 4 or less (see
\cite{babstein}).  Although all the possible Mahonian statistics based
on generalized patterns of length at most 3 were listed in
\cite{babstein}, proofs were not given that all of them were indeed
Mahonian.  These proofs were later supplied by Foata and Zeilberger
\cite{FZbast}, who solved some of the conjectures with bijections, but
others using the ``Umbral Transfer Matrix Method'' of Zeilberger
\cite{zeil-umbral}.  The remaining conjectures in \cite{babstein},
concerning a slight generalization of generalized patterns, were
proved bijectively by Foata and Randrianarivony \cite{foata-rand}.

\section{Avoidance (and occurrences) of generalized patterns of length 3}\label{sec-gp3}

The study of \emm avoidance, of generalized patterns---along the lines
of the work of Simion and Schmidt \cite{simion-schmidt} for the
classical patterns---was initiated by Claesson in the paper \emm
Generalized pattern avoidance, \cite{cla-gpa}, where the enumeration
was done for avoidance of any single pattern of length 3 with exactly
one dash.  These patterns fall into three equivalence classes with
respect to avoidance.  One of these classes, consisting of the
patterns \bxca, \bxac, \caxb\ and \acxb, has the same avoidance as
both of the classes of classical patterns of length 3, and thus has
avoidance enumerated by the Catalan numbers.  In fact, avoiding
\bxca\ is \emm equivalent, to avoiding \bxcxa, as shown in
\cite{cla-gpa}, and likewise for the other three patterns.  This is
obviously true also for the patterns $\ba$ and $\bxa$, since the only
permutation with no descent and no inversion is the increasing
permutation $123\dots n$.  It is shown by Hardarson \cite {hardarson}
that this can not occur for patterns of length greater than three,
that is, two different patterns of length more than three cannot be
avoided by the same permutations.  The other two classes of one-dash
patterns of length three, consisting of the patterns
%
%
\axbc\ and \axcb\ and their respective sets of equivalent patterns,
are enumerated by the Bell numbers, counting partitions of sets.  This
was proved bijectively in each case in \cite{cla-gpa}.

The fact that the Bell numbers count permutations avoiding \axbc\ (and
some other generalized patterns) has interesting implications.
Namely, this shows that the Stanley-Wilf conjecture (proved by Marcus
and Tardos \cite{marcus-tardos} in 2004) does not hold for some
generalized patterns.  This ex-conjecture says that for any classical
pattern $p$ the number of permutations of length $n$ avoiding $p$ is
bounded by $C^n$ for some constant~$C$.  This is easily seen to fail
for the Bell numbers, whose exponential generating function is
$e^{(e^x-1)}$.  Apart from this, nothing seems to be known about
growth rates of the number of permutations avoiding generalized
patterns.

That the Bell numbers count the avoidance of \axbc\ also implies the
falseness for generalized patterns of another conjecture, by Noonan
and Zeilberger \cite{noonan-zeilberger} (but first mentioned by Gessel
\cite{gessel}), that the number of permutations avoiding a classical
pattern is polynomially recursive, that is, satisfies a recursion
$$
P_0(n)f(n)= P_1(n)f(n-1)+P_2(n)f(n-2)+\cdots+P_k(n)f(n-k)
$$ where $k$ is a constant and each $P_i$ is a polynomial.  This
conjecture, however, is largely believed to be false even for
classical patterns, although a counterexample is still missing.

Claesson \cite{cla-gpa} also enumerated the avoidance of three classes
of pairs of generalized patterns of length 3 with one dash each.
These turned out to be equinumerous with non-overlapping set
partitions (counted by the Bessel numbers), involutions and Motzkin
paths, respectively.  Enumerative equivalences among the class of
patterns corresponding to non-overlapping partitions, together with
the connection to set partitions, naturally led to the definition of
\emm monotone, partitions in~\cite{cla-gpa}.

In \cite{cla-man-enum-bast}, Claesson and Mansour then completed the
enumeration of permutations avoiding any pair of generalized patterns
with one dash each.  They also conjectured enumerative results for
avoidance of any set of three or more such patterns.  These
conjectures were proved for sets of size three by Bernini, Ferrari and
Pinzani \cite{bernini-al-3-bast}, and by Bernini and Pergola
\cite{bernini-pergola} for the sizes 4, 5 and 6, the remaining sizes
being rather trivial.

Elizalde and Noy \cite{elizalde-noy} treated the {\em dashless
patterns} (which they call ``consecutive patterns''), that is,
patterns with no dashes, and gave generating functions enumerating
permutations according to the number $k$ of occurrences of a pattern.
This is a much stronger result than enumerating permutations \emm
avoiding, a pattern, which is just the special case $k=0$.  In
particular, they enumerated the avoidance of both Wilf classes of
dashless patterns of length 3, and gave differential equations
satisfied by the generating functions for three of the seven Wilf
classes of patterns of length 4 (see the next section).  As mentioned
before, their result in the special case of avoidance of the pattern
$123$, was obtained already in the book {\em Combinatorial
Enumeration} by Goulden and Jackson \cite[Exercise
5.2.17a]{goulden-jackson}.

For subsets of two or more dashless patterns of length 3, Kitaev
\cite{kit-multi} and Kitaev and Mansour
\cite{kit-mans-multi-gen,kit-mans-simul-gen} gave direct formulas for
almost every case, and recursive formulas for the few remaining ones.
Examples of the formulas thus obtained are $C_k+C_{k+1}$ where $C_k$
is the $k$-th Catalan number, the central binomial coefficients $\ch
2n,n,$, and the Entringer numbers, which also count certain
permutations starting with a decreasing sequence and then alternating
between ascents and descents.

Thus, the avoidance of any set of generalized patterns of length 3 has
been understood, in most cases in the sense of explicit formulas or
generating functions, or at least in terms of recursively defined
functions.

In \cite{cla-mans-count}, Claesson and Mansour found the number of
permutations with exactly one, two and three occurrences,
respectively, of the pattern $\bxca$.  They used the connection
between continued fractions and Motzkin paths (see Section
\ref{section-context}) that had been used in \cite{mad} to give a
continued fraction for a generating function counting occurrences of
$\bxca$, among other things.  Later, Parviainen \cite{parviainen}
showed how to use bicolored Motzkin paths to give a continued fraction
counting permutations of length $n$ according to the number of
occurrences of $\bxca$.  He gave an algorithm for finding an explicit
formula, and gave this explicit formula (always a rational function in
$n$ times a binomial coefficient of the form $\ch2n,n-a,$) for each
$n\le8$.  Finally, Corteel and Nadeau \cite{corteel-nadeau} found a
bijective proof of the fact, first proved in \cite{cla-mans-count},
that the number of permutations of length $n$ with exactly one
occurrence of $\bxca$ is $\ch2n,n-3,$.  They exploited the connection
between generalized patterns and \emm permutation tableaux,.  For more
about that connection, see Section \ref{section-context}.

\section{Patterns of length 4}\label{sec-gp4}

The classical patterns of length 4 fall into three Wilf classes,
represented by the patterns $1\dd2\dd3\dd4$, $1\dd3\dd4\dd2$ and
$1\dd3\dd2\dd4$.  The avoidance of the first two has been solved (the
first by Gessel \cite{gessel}, the second by B\'ona
\cite{bona-1324}), but $1\dd3\dd2\dd4$ still remains to be understood,
although some noteworthy progress was recently made by Albert et al.\
in~\cite{albert-et-al}.

For generalized patterns of length 4 (other than the classical ones),
the situation is much more complicated than for those of length 3, as
is to be expected. There are 48 symmetry classes, and computer
experiments show that there are at least 24 Wilf classes, but their
exact number does not seem to have been determined yet.

As mentioned above, for the dashless patterns of length 4, Elizalde
and Noy~\cite{elizalde-noy} gave differential equations satisfied by
the generating functions for the number of occurrences of three out of
seven Wilf classes, namely the classes containing 1234, 1243 and 1342,
respectively.  The remaining classes are represented by 2413, 2143,
1324 and 1423.  Note that a Wilf class is a class of patterns with the
same \emm avoidance,, whereas Elizalde and Noy proved that in each
Wilf class of dashless patterns of length four, all patterns have the
same \emm distribution,.  That is, the number of permutations of
length $n$ with $k$ occurrences of a pattern is the same for two
patterns in the same Wilf class in this case, but not in general.

Kitaev \cite[Theorem 13]{kit-pogp} found an expression for the
exponential generating function (EGF) for the avoidance of $\sigma\dd
k$, where $\sigma$ is any dashless pattern and $k$ is larger than all
the letters in $\sigma$, in terms of the EGFs for avoidance of the
pattern $\sigma$.  In particular, if $\sigma$ is any dashless pattern
of length 3, this, together with the results of Elizalde and Noy
\cite{elizalde-noy}, yields explicit formulas for the EGFs for the
avoidance of $\sigma\dd4$, where $\sigma$ is any dashless pattern of
length 3 (although one of these formulas involves the integral
$\int_0^x e^{-t^2}dt$).  Since there are precisely two Wilf classes
for dashless patterns of length 3, this gives the avoidance of two
Wilf classes of patterns of length 4.

Also, Callan \cite{callan-eigen} has given two recursive formulas for
the number of permutations avoiding $31\dd4\dd2$.

These seem to be all the enumerative results so far for patterns of
length 4.  For the non-classical patterns we thus have formulas for
the avoidance of six Wilf classes out of (at least)~24, whereas for
the classical patterns of length~4 there are formulas for two out of
three classes.  However, there is no reason to assume that the
non-classical patterns should be harder to deal with than the
classical ones, so explicit results for the avoidance of more such
patterns should not be considered out of reach.

\section{Generalized patterns appearing in other contexts}\label{section-context}

In Section \ref{section-intro} we mentioned the connection between
dashless patterns as valleys, peaks, etc., in permutations, and
Flajolet's \cite{flajolet} generation of Motzkin paths by means of
continued fractions.  Using results from Flajolet's paper
\cite{flajolet}, Clarke, Steingr\'imsson and Zeng
\cite[Corollary~11]{mad} found a continued fraction capturing, among
other things, the distribution of permutations according to the number
of occurrences of $\bxca$.  This was made explicit in
\cite[Corollary~23]{cla-mans-count}.  A polynomial formula for the
joint distribution of descents and $\bxca$ was conjectured by
Steingr\'imsson and Williams (unpublished), after Williams
\cite[Corollary~5.3]{williams} had shown that formula to count
permutations according to \emm weak excedances, and \emm alignments,
which Williams was studying in connection with so called \emm
permutation tableaux, (for definitions see
\cite{stein-williams,williams}).  This conjecture (part of a much
larger conjecture later proved in \cite{stein-williams}) was first
proved by Corteel \cite{corteel-crossings}.  The formula is as follows
(see Corollary~30 in \cite{stein-williams}):

The number of permutations of length $n$ with $k-1$ descents and $m$
occurrences of the pattern $\bxca$ is equal to the coefficient of
$q^m$ in
\begin{equation}
          q^{-k^2} \sum_{i=0}^{k-1} (-1)^i [k-i]^n q^{ki} \left( {n
       \choose i} q^{k-i} + {n \choose i-1}\right).
\end{equation}
Here, $[k-i]$ is the $q$-bracket defined by $[m]=(1+q+\dots+q^{m-1})$.
This is the only known polynomial formula for the entire distribution
of a pattern of length greater than 2.  The two cases of length 2
correspond to the Eulerian numbers, counting descents, and the
coefficients of the $q$-factorial
$$
[n]!=(1+q)(1+q+q^2)\cdots(1+q+q^2+\cdots+q^{n-1}),
$$ which count permutations according to the number of inversions.  A
descent is an occurrence of the pattern $\ba$ and an inversion is an
occurrence of $\bxa$.

Moreover, this connection between patterns and permutation tableaux
also led to the discovery, by Corteel \cite{corteel-crossings} (see
also \cite{corteel-williams-asep, corteel-williams-markov}), of a
connection between the permutation tableaux and the partially
asymmetric exclusion process (PASEP), an important model in
statistical mechanics.  In particular, the distribution of
permutations of length $n$ according to number of descents and number
of occurrences of the pattern $\bxca$ equals a probability
distribution studied for the PASEP.

In \cite{cks}, a bijection is given between the permutations of length
$n$ avoiding both $\bxdaxc$ and $\cxaxdxb$ on one hand, and so called
$\btrees$ on the other.  The \btrees\ are rooted plane trees with
certain labels on their vertices.  These trees were defined by
Jacquard and Schaeffer \cite{jacq-schaeffer}, who described a
bijection from rooted nonseparable planar maps to a set of labeled
plane trees including the \btrees. These trees represent, in a rather
transparent way, the recursive structure found by Brown and Tutte
\cite{brown-tutte} on planar maps.  As it turns out, the bijection
given in \cite{cks} simultaneously translates seven different
statistics on the permutations to corresponding statistics on the
\btrees.  In fact, the permutations avoiding $\bxdaxc$ and $\cxaxdxb$
seem to be more closely related structurally to the $\btrees$---and
thus to the planar maps involved---than the \emm two-stack sortable
permutations, that had previously been shown to be in bijection with
the planar maps in question (see Section \ref{section-disguise}).
Earlier, Dulucq, Gire and West \cite{dulucq-gire-west} constructed a
generating tree for the permutations avoiding $\bxdxaxc$ and $\cxadxb$
that they showed to be isomorphic to a generating tree for rooted
nonseparable planar maps. (Clearly, permutations avoiding these two
patterns are equinumerous with the permutations avoiding $\cxaxdxb$
and $\bxdaxc$, treated in \cite{cks}.)  However, instead of the
pattern $\cxadxb$ they used a so-called barred pattern, which we treat
in Section \ref{section-disguise}, and they only showed their
bijection, which is different from the one in \cite{cks}, to preserve
two different statistics, rather than the seven statistics in
\cite{cks}.

\section{Generalized patterns in disguise}\label{section-disguise}

As mentioned above, generalized patterns have occurred implicitly in
several places in the literature, even before the systematic study of
classical permutation patterns.  However, they have also appeared as
so-called \emm barred, patterns, which sometimes, but not always, turn
out to be equivalent (in terms of avoidance) to some generalized
patterns.  An example of a barred pattern is $4\dd\bar2\dd1\dd3$.  A
permutation $\pi$ is said to avoid this pattern if it avoids the
pattern $3\dd1\dd2$ (corresponding to the unbarred elements
$4\dd1\dd3$) \emm except, where that pattern is part of the pattern
$4\dd2\dd1\dd3$.

Gire \cite{gire} showed that the \emm Baxter permutations,, originally
defined in a very different way \cite{chung-al-baxter}, are those
avoiding the two barred patterns $4\dd1\dd\bar{3}\dd5\dd2$ and
$2\dd5\dd\bar{3}\dd1\dd4$.  It is easy to show that avoiding
$4\dd1\dd\bar{3}\dd5\dd2$ is equivalent to avoiding $3\dd14\dd2$ and
avoiding $2\dd5\dd\bar{3}\dd1\dd4$ is equivalent to avoiding
$2\dd41\dd3$.
Thus, the Baxter permutations are precisely those that avoid both
$3\dd14\dd2$ and $2\dd41\dd3$.  In fact, this was pointed out in 
Erik Ouchterlony's thesis \cite[p. 5]{ouchterlony-thesis}.

The barred pattern $4\dd1\dd\bar3\dd5\dd2$ also shows up in
\cite{dulucq-gire-west}, where Dulucq, Gire and West treated so called
\emm nonseparable permutations, (bijectively related to rooted
nonseparable planar maps), which they characterized by the avoidance
of $\bxdxaxc$ and $4\dd1\dd\bar3\dd5\dd2$, the latter one being
equivalent to $\cxadxb$ as mentioned in the previous section.

%
%

Also, in \cite{mireille-butler}, Bousquet-M\'elou and Butler deal with
the pattern $2\dd1\dd\bar3\dd5\dd4$, avoiding which is easily shown to
be equivalent to avoiding $\bxadxc$.  Permutations avoiding that
pattern and $\axcxbxd$ are called forest-like permutations in
\cite{mireille-butler}.  It is also mentioned there that avoiding
$2\dd1\dd\bar3\dd5\dd4$ (and thus $\bxadxc$) is equivalent to avoiding
$\bxaxdxc$ \emm with Bruhat condition $(1\leftrightarrow4)$, in the
%
%
terminology of Woo and Yong, who conjectured \cite{woo-yang} that a
Schubert variety is locally factorial if and only if its associated
permutation avoids these two patterns.  That conjecture was proved in
\cite{mireille-butler}.

Not all barred patterns can be expressed in terms of generalized
patterns, however.  For example, West \cite{west-thesis} showed that
two-stack sortable permutations are characterized by the simultaneous
avoidance of $\bxcxdxa$ and the barred pattern
$2\dd\bar5\dd3\dd4\dd1$, and it is easy to show that no generalized
pattern is avoided by the same permutations as those avoiding
$2\dd\bar5\dd3\dd4\dd1$.  It is also easy to show that there is no
pair of generalized patterns of length 4 that is avoided by the same
permutations as those that avoid $\bxcxdxa$ and
$2\dd\bar5\dd3\dd4\dd1$.  A consequence of this is that two-stack
sortability of a permutation cannot be characterized by the avoidance
of a set of generalized patterns.  This is because there are precisely
two permutations of length four that are not two-stack sortable,
namely 2341 and 3241, so such a set would have to contain two
generalized patterns whose underlying permutations (obtained by
disregarding the dashes in the patterns) were 2341 and 3241.  It is
easy to check, by computer, that no such pair will do the job.

Also Callan~\cite{callan-equi} has shown that the number of
permutations avoiding $31\dd4\dd2$ is the same as the number avoiding
the barred pattern $3\dd\bar5\dd2\dd4\dd1$, although these are not the
same permutations.  (As mentioned above, Callan \cite{callan-eigen}
has given two recursive formulas for this number.)

An obvious open problem here is to determine when avoiding a barred
pattern is equivalent to avoiding a generalized pattern.

\section{Asymptotics}\label{sec-asymp}

As mentioned before, the Stanley-Wilf conjecture, proved by Marcus and
Tardos \cite{marcus-tardos}, says that the number of permutations of
length $n$ avoiding any given classical pattern $p$ is bounded by
$C^n$ for some constant $C$ depending only on $p$.

In \cite{elizalde-asymptotic}, Elizalde studies asymptotics for the
number of permutations avoiding some generalized patterns and
concludes that, in contrast to the classical patterns, there probably
is ``a big range of possible asymptotic behaviors.''  Although much
work remains to be done here, and although it is not clear how varied
this behavior can be, the extremes are already known.  Namely, whereas
the number of permutations of length $n$ avoiding a classical pattern
is bounded by $c^n$ for a constant $c$, Elizalde \cite[Theorem
  4.1]{elizalde-asymptotic} shows that the number $\alpha_n(\sigma)$
of permutations of length $n$ avoiding a dashless (or \emm
consecutive,) pattern $\sigma$ of length at least three satisfies 
$$c^n n!< \alpha_n(\sigma)< d^n n!
$$ 
for some constants $c$ and $d$ (where, clearly, $0<c,d<1$).

An interesting open question related to this is which generalized
patterns, apart from the classical ones, satisfy the Stanley-Wilf
conjecture.  It has been pointed out by Hardarson
\cite{hardarson-personal} that a pattern containing a block with at
least two letters $a$ and $b$, with $a<b$, and a letter $x$ in some
other block, with $x<a$ or $b<x$, can not satisfy the Stanley-Wilf
conjecture, because there will be at least as many permutations
avoiding it as there are permutations that avoid $\axbc$, and the
number of such permutations is known not to satisfy the Stanley-Wilf
conjecture (these are the Bell numbers).  Thus, essentially the only
open cases left are the patterns $\bxcxda$, $\cxbxda$, $\cxdaxb$ and
$\bxdaxc$, which all have the same avoidance for $n\le10$.  If these
patterns turn out not to satisfy the Stanley-Wilf conjecture, that
would prove the tempting conjecture that a generalized pattern
satisfies the Stanley-Wilf conjecture if and only if it is avoided by
the same permutations as the classical pattern with the same
underlying permutation.  As mentioned in Section \ref{sec-gp3}, it
has been shown by Hardarson \cite{hardarson} that this happens only
for patterns of length 3 or less.

\section{Further generalizations}\label{sec-further}

Kitaev, in \cite{kit-pogp}, introduced a further generalization of
generalized patterns (GPs), namely the \emm partially ordered
generalized patterns, or POGPs.  These are GPs where some letters may
be incomparable in size.  An example of such a pattern is
$3\dd12\dd3$, an occurrence of which consists of four letters, the
middle two adjacent (and in increasing order) and the first and the
last both greater than the middle two, with no condition on the
relative sizes of the first and the last letter.  Avoiding the pattern
$3\dd12\dd3$ is equivalent to avoiding both $\cxabxd$ and $\dxabxc$.
Indeed, an \emm occurrence, of $3\dd12\dd3$ is equivalent to an
occurerence of either $\cxabxd$ or $\dxabxc$.  In general, a POGP is
equal, as a function counting occurrences, to a sum of GPs.

In \cite{hardarson}, Hardarson finds EGFs for the avoidance of
$k\dd\sigma\dd k$, where $\sigma$ is any dashless partially ordered
pattern and $k$ is larger than any letter in $\sigma$, in terms of the
EGF for the avoidance of $\sigma$.  In the special case where
$\sigma=12$ he gives a bijection between permutations avoiding
$3\dd12\dd3$ and bicolored set partitions, that is, all partitions of
a set where each block in the partition has one of two possible
colors. Also, for $\sigma=121$, he gives a bijection between
permutations of length $n+1$ avoiding $3\dd121\dd3$ and the Dowling
lattice on $\{1,2,\ldots,n\}$.

One interesting, and curious, result arising from Kitaev's study of
POGPs in~\cite{kit-pogp} is that knowing the EGF for the avoidance of
a dashless pattern $p$ is enough to find the EGF for the entire
distribution of the maximum number of \emm non-overlapping
occurrences, of $p$.  Two occurrences of $p$ in a permutation $\pi$
are non-overlapping if they have no letter of $\pi$ in common.  For
example, the permutation 4321 has three descents (occurrences of
$\ba$), but only two non-overlapping descents.  Namely, Kitaev
\cite{kit-pogp} proves the following theorem:
\begin{theorem*}[Kitaev \cite{kit-pogp}, Theorem 32]\label{kit-main}
Let $p$ be a dashless pattern. Let $A(x)$ be the EGF for the number
of permutations that avoid $p$ and let $N(\pi)$ be the
maximum number of non-overlapping occurrences of $p$ in $\pi$.
Then
$$
\displaystyle\sum_{\pi}y^{N(\pi)}\frac{x^{|\pi|}}{|\pi|!}=
\frac{A(x)}{1-y((x-1)A(x)+1)},
$$
where the sum is over all permutations of lengths $0,1,2,\dots$.
\end{theorem*}
Alternative proofs, and some extensions, of Theorem~\ref{kit-main}
were given by Mendes~\cite{mendes} and by Mendes and
Remmel~\cite{mendes-remmel}, using the theory of symmetric functions.

A generalization in a different direction is the study of generalized
patterns on words, that is, on permutations of multisets.  Research in
this area has only recently taken off, even in the case of classical
patterns.  For the generalized patterns, see \cite{bernini-al-words,
burs-mans-words, kit-gen-add}.

Finally, Burstein and Lankham have considered barred generalized
patterns, in relation to patience sorting problems (see
\cite{bur-lank-restricted}, which also has further references).

\end{document}